\documentstyle[11pt,epsf,amscd]{amsart}
\newtheorem{Proposition}{Proposition} [section]
\newtheorem{Theorem}{Theorem} [section]
\newtheorem{Corollary}{Corollary} [section]
\def\proof{\par{\it Proof}. \ignorespaces}
\def\endproof{{\ \vbox{\hrule\hbox{%
   \vrule height1.3ex\hskip0.8ex\vrule}\hrule }}\par}
\newenvironment{Proof}{\proof}{\endproof}

\begin{document}

\title{Topology of the iso-spectral real manifolds associated with
the generalized Toda lattices on semisimple Lie algebras
}

\author{Luis Casian$^*$ \ \
 and  \ Yuji Kodama$^\dagger$ \\ }

\thanks{*Department of Mathematics, The Ohio State University,
Columbus, OH 43210\endgraf
{\it E-mail address\/}: casian@@math.ohio-state.edu
\endgraf$\dagger$ Department of Mathematics, The Ohio State University,
Columbus, OH 43210\endgraf
{\it E-mail address\/}: kodama@@math.ohio-state.edu}

\begin{abstract}
This paper concerns the topology of isospectral real manifolds
of certain Jacobi elements associated with real split semisimple Lie
algebras.  The manifolds are related to the compactified level sets
of the generalized (nonperiodic) Toda lattice equations defined on
the semisimple Lie algebras.  We then give a cellular
decomposition and
the associated chain complex of the manifold by introducing colored
Dynkin diagrams which parametrize the cells in the decomposition.
We also discuss the Morse chain complex of the manifold.

\end{abstract}

\maketitle

\markboth{LUIS CASIAN AND YUJI KODAMA}
  {TOPOLOGY OF TODA LATTICES}

\section{The generalized Toda lattice equations}
\renewcommand{\theequation}{1.\arabic{equation}}\setcounter{equation}{0}
\renewcommand{\thefigure}{1.\arabic{figure}}\setcounter{figure}{0}

Let $\frak g$ denote a real
split semisimple Lie algebra of rank $l$.  We fix a
split Cartan subalgebra $\frak h$ with root system $\Delta $, real root vectors
$e_{\alpha_i}$ associated with simple roots $\{ \alpha_i : i=1,..,l \}=\Pi$.
We also denote $\{h_{\alpha_i},e_{\pm\alpha_i}\}$ the Cartan-Chevalley basis of
$\frak g$ which satisfies the relations,
\begin{equation}
\label{chevalley}
 [h_{\alpha_i} , h_{\alpha_j}] = 0, \, \,
  [h_{\alpha_i}, e_{\pm \alpha_j}] = \pm C_{j,i}e_{\pm \alpha_j} \ , \, \,
  [e_{\alpha_i} , e_{-\alpha_j}] = \delta_{i,j}h_{\alpha_j}.
\end{equation}
where the $l\times l$ matrix
$(C_{i,j})$ is the Cartan matrix
corresponding to $\frak g$, and $C_{i,j}=\alpha_i(h_{\alpha_j})=
\langle \alpha_i,h_{\alpha_j}\rangle$.

 Then the generalized Toda lattice equation related to
real split semisimple Lie algebra is defined by
the following system of 2nd order differential equations for
the real variables $\{f_j(t) : j=1,\cdots,l\}$,
\begin{equation}
\label{toda}
 {d^2 f_i \over dt^2}= \epsilon_i \exp \left(
-\langle \alpha_i, f\rangle \right),
\end{equation}
where $f=\sum_{j=1}^l f_j(t)h_{\alpha_j} \in {\frak h}$ and
$\epsilon_i\in \{\pm 1\}$.

\vskip 0.5cm
\noindent
{\bf Remark 1.1.}
The case with ${\frak g}={\frak {sl}}(l+1,{\Bbb R})$ corresponds
to the indefinite Toda lattice introduced in \cite{kodama:96}.
The main feature of the indefinite Toda equation having
at least one of $\epsilon_i$ being $-1$ is that the
solution blows up to infinity in finite time \cite{kodama:96}.
Having introduced the signs, the group corresponding to the Toda lattice is
a real split Lie group $\tilde G$ with Lie algebra
$\frak g$.
 For example, in the case of ${\frak g}={\frak{sl}}(n,\Bbb R)$,
if $n$ is odd, $\tilde
G=SL(n,\Bbb R)$, and if $n$ is even, $\tilde G=Ad( SL(n,\Bbb R)^{\pm})$.

\vskip 0.5cm
\noindent
{\bf Remark 1.2.}
If we consider the complex Toda equation, $\epsilon_i$ in
(\ref{toda}) can be absorbed in $f_i\in {\Bbb C}$, so that the
present study deals with the diconnected Cartan subgroup,
where the generalized Toda lattice defines a flow in
each connected component.

\vskip 0.5cm
\noindent
{\bf Remark 1.3.}
The original Toda lattice in \cite{toda:67}
is obtained as the case with all $\epsilon_i=1$
where the position of the $i$-th particle is given by
$q_i=f_i-f_{i+1}$ for $i=1,\ldots,l$,
\begin{equation}
\label{otoda}
{d^2 q_i \over dt^2}=  \exp (q_{i-1}-q_{i})-\exp (q_i-q_{i+1}),
\end{equation}
where $f_{l+1}=0$ and $f_{0}=f_{l+2}=-\infty$
indicating $q_0=-\infty$ and $q_{l+1}=\infty$.

\subsection{Lax formulation: isospectral manifold $Z(\gamma)_{\Bbb R}$}

The system (\ref{toda}) can be written in a Lax equation
which describes an iso-spectral deformation of a Jacobi element
of $\frak g$ \cite{flaschka:74}.
Define the set of real functions $\{(a_i(t), b_i(t)) :
i=1,\ldots,l\}$,
\begin{equation}
 \displaystyle{a_i(t)={d \over dt}f_i(t)},~~~~~~
 \displaystyle{b_i(t)=\epsilon_i \exp \left(-\langle \alpha_i,
f\rangle\right)}.
\label{ab}
\end{equation}
Then the Toda equation (\ref{toda}) can be written in the Lax form
\cite{flaschka:74,kostant:79},
\begin{equation}
\label{lax}
\displaystyle{{dX \over dt}=[P, X]}
\end{equation}
where the Lax pair $(X,P)$ are defined by
\begin{equation}
\left\{
\begin{array} {ll}
& \displaystyle{X(t)=\sum_{i=1}^l a_i(t)h_{\alpha_i}+\sum_{i=1}^l \left(
b_i(t)e_{-\alpha_i}+e_{\alpha_i}\right)} \\
& \displaystyle{P(t)=-\sum_{i=1}^l b_i(t)e_{-\alpha_i}}
\end{array}
\right.
\label{xp}
\end{equation}
The Lax form (\ref{lax}) represents an isospectral deformation of
the Jacobi element $X$.

We denote the disconnected manifold given by the set of the elements
in the form $X$ of $\frak g$,
\begin{equation}
Z_{\Bbb R}=\left\{X=x+\sum_{i=1}^l(e_{\alpha_i}+b_i e_{-\alpha_i})
\in{\frak g}~
:~ x\in {\frak h},~ b_i\in {\Bbb R}^*\right\}
= \bigcup_{\epsilon\in {\cal E}} Z_{\epsilon},
\label{zr}
\end{equation}
where ${\Bbb R}^*={\Bbb R}\setminus \{0\}$,
and the connected component $Z_{\epsilon}$ is given by
\begin{equation}
\label{zepsilon}
Z_{\epsilon}=\{ X\in Z_{\Bbb R}~ :~ \epsilon=(\epsilon_1,\ldots,\epsilon_l)
\in {\cal E},
~ sign(b_i)=\epsilon_i\}.
\end{equation}
Here ${\cal E}$ is the set of all the signs $\epsilon$, so the set $Z_{\Bbb R}$
is the disjoint union of the $2^l$ connected components.

A real isospectral leaf in
$Z_{\Bbb R}$ is defined by the level sets of the
Chevalley invariants, denoted as $(I_1,\ldots,I_l)$,
which are the polynomials of the variables $(a_i,b_i)$.
The invariants then define a differentiable map,
\begin{eqnarray}
\nonumber
{\cal I} &:& Z_{\Bbb R} \longrightarrow {\Bbb R}^l\\
        &{}& X \longmapsto \gamma=(I_1,\ldots,I_l)
\label{I}
\end{eqnarray}
The real isospectral leaf $Z(\gamma)_{\Bbb R}$ is then given by
\begin{equation}
\label{Z}
Z(\gamma)_{\Bbb R}={\cal I}^{-1}(\gamma)\bigcap Z_{\Bbb R},
\end{equation}
Our main purpose in this paper is to give a detailed structure of
the compactified manifold $\hat Z(\gamma)_{\Bbb R}$ from
the viewpoint of the Lie group theory.

\vskip 0.5cm
\noindent
{\bf Remark 1.4.}
The Chevalley invariants provide $l$-involutive integrals for
the generalized Toda lattice equation, so that this proves
the integrability of the equation in the Liouville-Arnold
sense.

\begin{figure}
\label{fig:1}
\epsfysize=8cm
\centerline{\epsffile{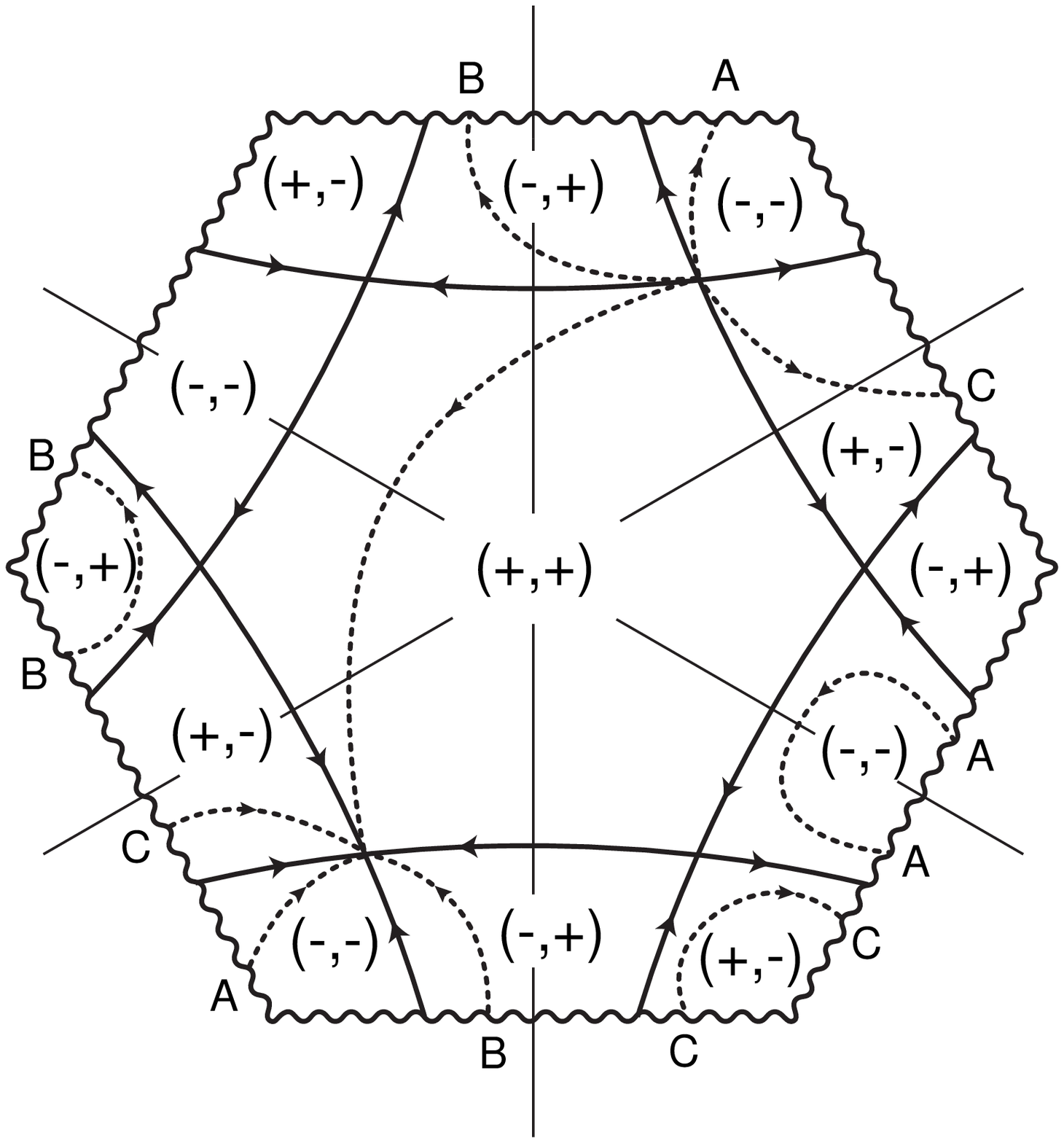}}
\caption{The isospectral manifold ${\hat Z}(\gamma)_{\Bbb R}$ for
${\frak{sl}}(3, {\Bbb R})$.}
\end{figure}

\vskip 0.5cm
\noindent
{\bf Remark 1.5.}
The construction of the compactified manifold $\hat Z(\gamma)_{\Bbb R}$
for the case of ${\frak{sl}}(l+1,{\Bbb R})$ was given in
\cite{kodama:98} based on the explicit solution structure in
terms of the $\tau$-functions, which provide a local coordinate
system for the manifold. By tracing the solution orbit
of the indefinite Toda equation, the disconnected components in
$Z(\gamma)_{\Bbb R}$ are all
glued together to make a smooth compact manifold.
 The result is maybe
well explained by Figure \ref{fig:1} for the case of
$A_2\cong{\frak{sl}}(3,{\Bbb R})$.
In the figure, the Toda orbits are shown as the dotted lines, and
each region labeled by the same signs in
$(\epsilon_1,\epsilon_2)$ with $\epsilon_i \in
\{\pm \}$ are glued together through the boundary (the wavy-lines) of
the hexagon. At a point of the boundary the Toda orbit
blows up in finite time, but the orbit can be uniquely
traced to the one in the next region (marked by the same letter
$A, B$ or $C$). Also the flows on the solid lines show
the solutions of the subsystems (i.e. either $b_1=0$ or $b_2=0$).
Then the compactified manifold $\hat Z(\gamma)_{\Bbb R}$
by adding the blow-up points (the wavy-lines)
 and the flows of the subsystems (the solid-lines)
to $Z(\gamma)_{\Bbb R}$
in this case is shown to be isomorphic to
the connected sum of two Klein bottles, that is, the integral
homology $H_k({\hat Z}(\gamma)_{\Bbb R}, {\Bbb Z})$  is given by
$H_0={\Bbb Z}, H_1={\Bbb Z}^3\oplus {\Bbb Z}_2$, and $H_2=0$.
In the case of
${\frak{sl}}(n, {\Bbb R})$ for $n\ge 3$, $\hat Z(\gamma)_{\Bbb R}$ is shown
to be
nonorientable and the symmetry group is the semi-direct product
of $({\Bbb Z}_2)^{n-1}$ and the Weyl group $W=S_n$, the permutation group.
One should compare this
with the result of Tomei \cite{tomei:84} where the compact manifolds are
associated with
the definite (original) Toda lattice equation and the compactification
is done by adding only the subsystems. (Also see \cite{davis:98} for some
topological aspects of the manifolds which are identified as permutohedrons.)

\subsection{Leznov-Saveliev formulation: Cartan subgroup $H_{\Bbb R}$}

In the zero curvature formulation in \cite{leznov:79},
the generalized Toda lattice equation (\ref{toda}) with the sign
$\epsilon=(\epsilon_1,\ldots,\epsilon_l)$
can be expressed as an orbit on the connected component $H_{\epsilon}$
of the Cartan subgroup $H_{\Bbb R}$,
\begin{equation}
\label{H}
H_{\Bbb R}=\bigcup_{\epsilon\in {\cal E}} H_{\epsilon}
\end{equation}
where $H_{(1,\ldots,1)}:=H=\exp{\frak h}$, the connected component
with the identity.  Thus the set $H_{\Bbb R}$ consists of
$2^l$ connected components.  Let $g_{\epsilon}$ be an element
of $H_{\epsilon}$ given by
\begin{equation}
\label{ge}
g_{\epsilon}=h_{\epsilon}\exp f
\end{equation}
which can be also considered as a map
from $Z_{\Bbb R}$ to $H_{\Bbb R}$.  Here the element
 $h_{\epsilon}\in H_{\epsilon}$ satisfies
$\chi_{\alpha_i}(h_{\epsilon})=\epsilon_i$ with the group
character $\chi_{\phi}$ determined by a root $\phi\in \Delta$,
and each connected component of $H_{\Bbb R}$ can be written as
$H_{\epsilon}=h_{\epsilon}H$.
Then the Toda lattice (\ref{toda}) is written as an evolution of
$g_{\epsilon}(t)$,
\begin{equation}
\label{lez}
{d \over dt}g_{\epsilon}^{-1}{d \over dt}g_{\epsilon}
=\left[g_{\epsilon}^{-1}e_{+} g_{\epsilon} ~, ~e_{-}\right],
\end{equation}
where $e_{\pm}$ are fixed elements in the simple root spaces
${\frak g}_{\pm\Pi}$ so that all the elements in ${\frak g}_{\pm\Pi}$
can be generated by $e_{\pm}$, i.e.
${\frak g}_{\pm\Pi}=\{Ad_h(e_{\pm}): h\in H\}$.  In particular,
we take
\begin{equation}
\label{e}
e_{\pm}=\displaystyle{\sum_{i=1}^n e_{\pm\alpha_i}}.
\end{equation}
With the group character $\chi_{\alpha_i}$,
the solution $b_i(t)$ of the Toda lattice is given by
\begin{equation}
\label{bi}
b_i(t)=\left[\chi_{\alpha_i}(g_{\epsilon})\right]^{-1}=\chi_{-\alpha_i}
(g_{\epsilon}).
\end{equation}

\vskip 0.5cm
\noindent
{\bf Remark 1.6.}
With the fundamental weights $\omega_i$ defined as
$\langle \omega_i, h_{\alpha_j}\rangle
=\delta_{ij}$, i.e. $\alpha_i=\sum_{j=1}^lC_{i,j}\omega_j$,  we can write
the solution
\begin{equation}
\label{tau}
b_i=\prod_{j=1}^l\left[\chi_{\omega_j}(g_{\epsilon})\right]^{-C_{i,j}},
\end{equation}
which is the well-known $\tau$-function representation of the
solution with $\tau_i(t):=\chi_{\omega_i}(g_{\epsilon})$.

\vskip 0.5cm
\noindent
{\bf Remark 1.7.}
In the compactification of the disconnected Cartan subgroup $H_{\Bbb R}$,
we need to add pieces corresponding to the
blow-ups ($|b_i|=\infty$) and the subsystems ($b_i=0$).
The subsystems are determined by the subset $A=\{\alpha_i\in\Pi:
b_i\ne 0\}$, and the
corresponding
Cartan subgroup, denoted by $H^A_{\Bbb R}$, may be defined as
\begin{equation}
\label{hsub}
H_{\Bbb R}^A=\bigcup_{\epsilon \in {\cal E}^A} h_{\epsilon}H^A,
\end{equation}
where the set ${\cal E}^A\subset {\cal E}$ and the Cartan subgroup $H^A$
are defined by
\begin{equation}
{\cal E}^A=\{(\epsilon_1,\ldots,\epsilon_l)\in {\cal E}~:~
\epsilon_i=1~{\mbox{if}}~ \alpha_i\in A \}
\end{equation}
\begin{equation}
H^A=\exp {\frak h}^A~~\quad{\mbox{with}}~\quad{\frak h}^A
={\mbox{Span}}_{\Bbb R}\{h_{\alpha_i}\in {\frak h}:\alpha_i\notin A\}.
\end{equation}
Then the subsystems are also expressed as the same form of (\ref{lez}) with
$g_{\epsilon}^A\in H_{\epsilon}^A$,
\begin{equation}
\label{gea}
g^A_{\epsilon}=h_{\epsilon}\exp\left(\sum_{\alpha_i\notin A}
f_i(t)h_{\alpha_i}\right).
\end{equation}
The corresponding Lax pair ($X^A, P^A$) is given by
\begin{equation}
\label{xpa}
\left\{
\begin{array}{ll}
&\displaystyle{ X^A=\sum_{\alpha_i\in \Pi}a_ih_{\alpha_i}+
\sum_{\alpha_i\notin A} b_ie_{-\alpha_i} + e_+} , \\
&\displaystyle{P^A=-\sum_{\alpha_i\notin A} b_ie_{-\alpha_i}}.
\end{array}
\right.
\end{equation}
which is just the Lax pair (\ref{xp}) with $b_{i}=0$ for $\alpha_i\in A$.
Note here that $a_i(t)=$constant if $b_i=0$.
We also consider that the dimension of the manifold $H_{\Bbb R}^A$
is $l-|A|$, the number of parameters $f_i$.

\vskip 0.5cm
\noindent
{\bf Remark 1.8}
The compactification of the isospectral manifold $Z(\gamma)_{\Bbb R}$ for a
fixed
$\gamma \in {\Bbb R}^l$ can be obtained by sending it  to the flag manifold
${\tilde G}/B_+$ with the Borel subgroup $B_+$ of $\tilde G$ \cite{kostant:79},
so that the compactified manifold $\hat Z(\gamma)_{\Bbb R}$ is a toric
variety $\overline{H_{\Bbb R}x B_+}$ with a generic element $x\in \tilde G$.
Then we can show:
\begin{Theorem}
The isospectral manifold $\hat Z(\gamma)_{\Bbb R}$ is a smooth compact manifold
diffeomorphic to $\hat H_{\Bbb R}$.
\end{Theorem}
The complex version of this theorem is given in \cite{flaschka:91},
and the proof of the present case is essentially given in the same manner
(the detail of the proof is given in \cite{casian:99}).

In the following two sections, we will describe the structure of
$\hat H_{\Bbb R}$ using the Weyl action on the manifold. This is a brief
summary of a preprint \cite{casian:99}, and the proofs of the
results (Proposition 2.1, Theorem 3.1 and 3.2) can be found there.
Then in the section 4, we will present the Morse theory to compute
the integral homology of the manifold $\hat H_{\Bbb R}$.

\section{The structure of $\hat H_{\Bbb R}$ as the union of the subsystems}
\renewcommand{\theequation}{2.\arabic{equation}}\setcounter{equation}{0}

As was shown in the previous section, the set $H_{\Bbb R}$ can be
parametrized by the group characters $\chi_{\alpha_i}$, that is,
$H_{\Bbb R}=\cup_{\epsilon\in {\cal E}} H_{\epsilon}$ with $H_{\epsilon}=
h_{\epsilon}H$,
\begin{equation}
\label{hepsilon}
H_{\epsilon}=\left\{\ h\in H_{\Bbb R}~:~sign(\chi_{\alpha_i}(h))=
\chi_{\alpha_i}(h_{\epsilon})=\epsilon_i,
~{\mbox{for}}~ i=1,\ldots,l\right\}.
\end{equation}
Note that each $H_{\epsilon}$ is diffeomorphic to ${\Bbb R}^l$.

Since the Weyl group acts on $H_{\Bbb R}$, one can partition $H_{\Bbb R}$
into the $|W|$ convex cones of the Weyl chambers. We denote the cone
in the antidominant chamber as
\begin{equation}
\label{dominant}
H_{\Bbb R}^-=\bigcup_{\epsilon\in{\cal E}}H^-_{\epsilon}
\end{equation}
where the connected component $H_{\epsilon}^-$ is defined by
\begin{equation}
\label{he-}
H_{\epsilon}^-=\left\{h \in H_{\Bbb R}~:~ |\chi_{-\alpha_i}(h)|\le 1,
~ sign(\chi_{\alpha_i}(h))=\epsilon_i \right\}
\end{equation}
The boundaries of he chamber $H_{\epsilon}^-$ corresponding to
$\chi_{\alpha_i}(h)=1$ and $\chi_{\alpha_i}(h)=-1$ are called
the positive and negative $\alpha_i$-walls, and especially
the positive $\alpha_i$-wall gives the hyperplane of
the Weyl-reflection with respect to the root $\alpha_i$.
Then the connected component $H_{\epsilon}$ of $H_{\Bbb R}$
is expressed as the union of $W$-translations of $H_{\epsilon}^-$,
i.e.
\begin{equation}
\label{w-he}
H_{\epsilon}=\bigcup_{w\in W}w\left(H^-_{\epsilon(w)}\right).
\end{equation}
Here the $W$-action on $H_{\epsilon}^-$ is obtained through
the action on the group characters,
\begin{equation}
\label{W}
s_{\alpha_i}(\chi_{\alpha_j})=\chi_{s_{\alpha_i}\alpha_j}=\chi_{\alpha_j}
\chi_{\alpha_i}^{-C_{j,i}},
\end{equation}
from which we also have the $W$-action on the set ${\cal E}$ as
$ s_{\alpha_i}: \epsilon_j \mapsto \epsilon_j'$ with $\epsilon_j
=\chi_{\alpha_j}(h_{\epsilon})$,
\begin{equation}
\label{ewj}
\epsilon_j'=\epsilon_j\epsilon_i^{-C{j,i}},
\end{equation}
where we have used the Weyl-reflection,
$s_{\alpha_i}\alpha_j=\alpha_j-C_{j,i}\alpha_i$.
Thus the element $h_{\epsilon}$ is $W$-translated to $h_{\epsilon'}$
with the sign $\epsilon'=(\epsilon_1',\ldots,\epsilon_l')$
given by (\ref{ewj}), which we denote $\epsilon'=\epsilon(w)$ with
$w=s_{\alpha_i}$-Weyl reflection.  Thus, with the decomposition
(\ref{w-he}), we can consider only the antidominant chamber $H_{\epsilon}^-$,
and obtain the whole $H_{\Bbb R}$ by the $W$-translates.
This is also true for the compactified manifold $\hat H_{\Bbb R}$.

Let us first make the closure of $H_{\epsilon}^-$ by adding the pieces
corresponding
to the subsystems having the lower dimensions $l-|A|$
where $A\subset \Pi$ determines the subsystem (see Remark 1.7).
We let
\begin{equation}
\label{hea-}
H_{\epsilon}^{A,-}=\left\{h\in H_{\epsilon}^A~:~|\chi_{-\alpha_i}(h)|\le 1,
~\epsilon=(\epsilon_{i_1},\ldots,\epsilon_{i_m})~{\mbox{for}}~
\alpha_{i_j}\notin A\right\}
\end{equation}
Then the closure of the set $H_{\epsilon}^-$ can be obtained by
\begin{equation}
\label{he-c}
{\overline{H_{\epsilon}^-}}=\bigcup_{A\subset \Pi} H_{\epsilon}^{A,-},
\end{equation}
and the compactified manifold $\hat H_{\Bbb R}$ is given by the
$W$-translates of (\ref{he-c}), i.e.
\begin{equation}
\label{hrhat}
{\hat H}_{\Bbb R}=\bigcup_{\epsilon\in{\cal E}}
\bigcup_{w\in W}w\left({\overline{H_{\epsilon(w)}^-}}\right).
\end{equation}
We summarize the result as:
\begin{Proposition}
The closed set $\overline{H_{\epsilon}^-}$ is isomorphic to the box
$\{(t_1,\ldots,t_l)\in {\Bbb R}^l: -1\le t_j\le 1\}$, and
the manifold $\hat H_{\Bbb R}$ is compact and has an action of
the Weyl group $W$.
\end{Proposition}

\section{Topology of $\hat H_{\Bbb R}$}
\renewcommand{\theequation}{3.\arabic{equation}}\setcounter{equation}{0}
\renewcommand{\thefigure}{3.\arabic{figure}}\setcounter{figure}{0}

\subsection{Colored Dynkin diagrams}

We here give a cellular decomposition and construct the associated
chain complex of the compactified manifold $\hat H_{\Bbb R}$.
We first introduce the set of colored Dynkin diagrams to parametrize
the cells in the decomposition.  A colored Dynkin diagram is simply a
Dynkin diagram
in which some of the vertices have been colored either red ($R$) or blue ($B$).
For examples, in the case of $A_2\cong {\frak{sl}}(3,\Bbb R)$,
we have $\circ_R-\circ$, $\circ_B-\circ_R$, etc.
Thus a colored Dynkin diagram $D$ corresponds to a pair $(S, {\eta})$
with $S\subset\Pi$ and $\eta:S\to \{\pm1\}$, where $\eta(\alpha_i)=-1$
if $\alpha_i$ is colored $R$, and $\eta(\alpha_i)=1$ if $\alpha_i$
is colored $B$.  We denote the set of colored Dynkin diagrams as
\begin{equation}
\label{cd}
{\Bbb D}(S)=\left\{ D=(S,\eta) ~ :~ S\subset \Pi, ~\eta(\alpha)\in \{\pm 1\}~
{\mbox{for}}~ \alpha\in S\right\}.
\end{equation}

Let $W_S$ be the group generated by the simple reflections corresponding
to the roots in $S$.
We then define the $W_S$-action on the set ${\Bbb D}(S)$ as follows:
For any $\alpha_i\in S$, $s_{\alpha_i}D=D'$ is a new colored Dynkin diagram
having the colors corresponding to the sign change $\epsilon_j'=\epsilon_j
\epsilon_i^{-C_{j,i}}$ in (\ref{ewj}) with the identification that $R$
if the sign is $-1$, and $B$ if it is $+1$.  For example, in the case of $A_2$,
we have $s_{\alpha_1}(\circ_R-\circ_B)=\circ_R-\circ_R,
s_{\alpha_1}(\circ_B-\circ_R)=\circ_B-\circ_R$.

The $W_S$-action induces the $W$-translates on the set
${\Bbb D}(S)$ as $W\times_{W_S}{\Bbb D}(S)$.  The elements of this set
are given by pairs $(w, D)$, and among the elements we have an
equivalence relation $\sim$, that is, $(wx, D)\sim (w, xD)$ for
any $x\in W_S$.  The equivalent relation gives a bijective correspondence
between $W\times_{W_S}{\Bbb D}(S)$ and ${\Bbb D}(S)\times W/W_S$.
We then define the set ${\Bbb D}^k$ as
\begin{equation}
\label{dk}
{\Bbb D}^k:=\left\{ (D, [w]_{\Pi-S})~:~ D\in {\Bbb D}(S), ~
[w]_{\Pi-S}\in W/W_S, ~ |S|=k\right\},
\end{equation}
which parametrizes all the connected components of the Cartan subgroups of
the form $H_{\Bbb R}^{\Pi-S}$ corresponding to the subsystems
defined in Remark 1.7.
In this parametrization, ${\Bbb D}^k$ corresponds explicitely to
the dual of the set $H_{\Bbb R}^{\Pi-S}$, so that the parametrized cell
has the codimension $k$, and dim$H_{\Bbb R}^{\Pi-S}=k$.
Thus all the cells in $\hat H_{\Bbb R}$ can be parametrized by the
sets ${\Bbb D}^k$, and we have:
\begin{Theorem}
The collection of the sets ${\Bbb D}^k$ defined (\ref{dk}) gives a cell
decomposition of the compact manifold $\hat H_{\Bbb R}$.
\end{Theorem}

\vskip 0.5cm
\noindent
{\bf Remark 3.1}
There is a more convenient cell decomposition of $\hat H_{\Bbb R}$
for the purpose of calculating homology explicitly. The only change is that
the $l$ dimensional cell becomes the union of all the $l$-cells together with
all the (internal) boundaries corresponding to colored Dynkin diagrams
where all the colored vertices are colored $B$. This is the set,
$$\hat H_{\Bbb R} \setminus
\bigcup_{S\subset\Pi,~ w\in W \atop \eta(\alpha_i)=-1
 ~for~some~\alpha_i\in S}
 (S,\eta,[w]_{\Pi-S}).$$
This set can be seen to be homeomorphic to
 ${\Bbb R}^l$. With this cell decomposition there is exactly one  $l$-cell,
 and the other lower dimensional cells correspond to colored Dynkin diagrams
 $D$, in which at least one vertex of $D$ has been colored $R$.

\begin{figure}
\label{fig:2}
\epsfysize=8cm
\centerline{\epsffile{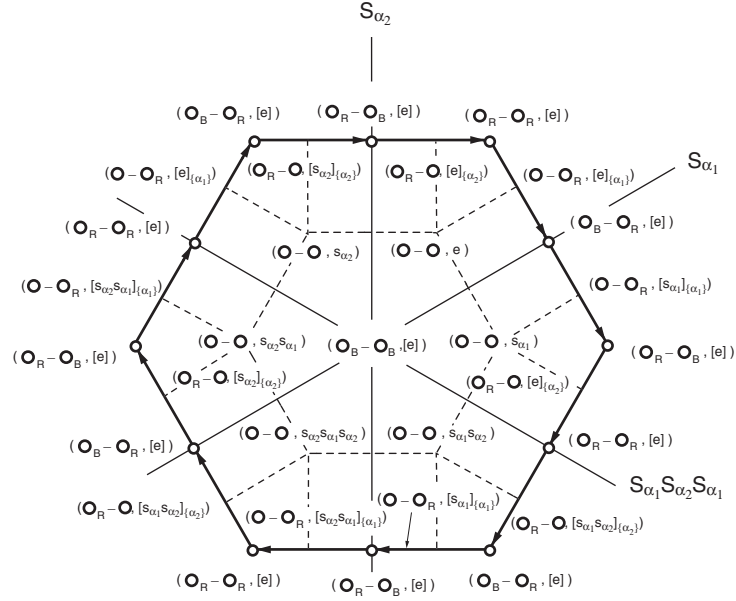}}
\caption{The manifold ${\hat H}_{\Bbb R}$ parametrized by colored Dynkin
diagrams
for ${\frak{sl}}(3, {\Bbb R})$.}
\end{figure}

\vskip 0.5cm
\noindent
{\bf Example 3.1.}
In the case of $A_2$, we have:
\begin{enumerate}
\item{ For $k=2$, i.e. $S=\Pi$, we have 4 vertices (0-cell) parametrized by
the elements of ${\Bbb D}^2$, which correspond to the 4 connected components
of $H_{\Bbb R}$ as dual cells,
$$
(\circ_B-\circ_B, [e]),~(\circ_B-\circ_R, [e]),~(\circ_R-\circ_B, [e]),
~(\circ_R-\circ_R, [e]),~
$$
where $[e]=[e]_{\emptyset}=[w]_{\emptyset}$ for any $w\in W$.}
\item{ For $k=1$, if $S=\{\alpha_1\}$, we have 6 1-cells
parametrized by,
$$
(\circ_B-\circ, [w]_{\{\alpha_2\}}),~(\circ_R-\circ, [w]_{\{\alpha_2\}}),~
$$
where $W/W_S=\{e, s_{\alpha_2}, s_{\alpha_1}s_{\alpha_2}\}$, and if
$S=\{\alpha_2\}$,
we have also 6 1-cells,
$$
(\circ-\circ_B, [w]_{\{\alpha_1\}}),~(\circ-\circ_R, [w]_{\{\alpha_1\}}),~
$$
where $W/W_S=\{e,s_{\alpha_1}, s_{\alpha_2}s_{\alpha_1}\}$.
Those colored Dynkin diagrams with $w=e$ correspond to the
4 walls (2-posive and 2-negative walls) of the antidominant chamber
$H_{\Bbb R}^-$,
which is isomorphic to a {\it square}.}
\item{ For $k=0$, i.e. $S=\emptyset$, we have $6=|W|$ 2-cells of the convex
cones
corresponding to the Weyl chambers parametrized by
$$
(\circ-\circ, w)
$$ for $w\in W$.  Those are dual to the 6-vertices corresponding to
$H_{\Bbb R}^{\Pi}$. As mentioned in Remark 3.1, we have a simpler cell
decomposition.
Namely, the union of all colored Dynkin
diagrams having no $R$-colored vertices forms the unique 2-cell which is the
set of internal points of the hexagon homeomorphic to ${\Bbb R}^l$,
and all other cells consists of the
boundary of the hexagon.}
\end{enumerate}
Figure \ref{fig:2} illustrates the example. One should note that the full
parametrization
of the cells are obtained by the $W$-translates of $(D, [e]_{\Pi-S})$
corresponding to
the subsystems in the antidominant chamber $H_{\Bbb R}^{\Pi-S,-}$ for
all the choices of $S\subset \Pi$.

\subsection{Boundary maps and chain complex}

With those parametrizations of the cells, we now define the boundary maps
on the colored Dynkin diagrams.  Let $(j,c)$ be a pair of integers with
$j=1,\ldots,m$ and $c=1,2$.  Then we define the $(j,c)$-boundary,
denoted as $\partial_{j,c} D$, of
a colored Dynkin diagram $D$ with the set of unclored vertices
$\{\alpha_{i_j}:1\le i_1.\ldots,i_m\le l\}$ as
a new colored Dynkin diagram by coloring the $i_j$-th vertex with $R$
if $c=1$ and with $B$ if $c=2$.  Recall that a colored Dynkin diagram $D$
corresponds to a pair $(S,\eta)$ with $S\subset \Pi$ and $\eta:S\to \{\pm 1\}$.
Thus the boundary operator $\partial_{j,c}$ determines a new pair $(S\cup
\{\alpha_{i_j}\}, \eta')$ where $\eta'$ is an extension of the $\eta$ on
$S\cup\{\alpha_{i_j}\}$ with $\eta'(\alpha_{i_j})=(-1)^c$. Giving an
orientation
on the boundary, we define the map $\tilde \partial_{j,c}:=(-1)^{j+c+1}
\partial_{j,c}$ which gives a map on the ${\Bbb Z}$-modules,
\begin{equation}
\label{boundary}
\tilde \partial_{j,c}~:~ {\Bbb Z}[{\Bbb D}(S)] ~ \longrightarrow ~
{\Bbb Z}[{\Bbb D}(S\cup\{\alpha_{i_j}\})].
\end{equation}

Let us now consider the ${\Bbb Z}$-modules of the full set of colored Dynkin
diagrams, ${\Bbb D}^k={\Bbb D}(S)\times W/W_S |_{|S|=k}$.
We denote the module as
\begin{equation}
\label{wmod}
{\cal M}(S)={\Bbb D}[W]\otimes_{{\Bbb Z}[W_S] }{\Bbb Z}[{\Bbb D}(S)],
\end{equation}
so that the $k$-chain is given by the direct sum of all these
modules over all sets with $|S|=l-k$,
\begin{equation}
\label{mk}
{\cal M}_{k}=\bigoplus_{|S|=l-k}{\cal M}(S).
\end{equation}
The boundary map $\partial_k : {\cal M}_{k} \to {\cal M}_{k-1}$ can be
defined by
\begin{equation}
\label{bmod}
\partial_{k}\left(D,[w]_{\Pi-S}\right)=\sum_{1\le j\le k \atop{c=1,2}}
\left(\tilde \partial_{j,c} D, [w]_{\Pi-\{S\cup\{\alpha_{i_j}\}\}}\right),
\end{equation}
where $\Pi-S=\{\alpha_{i_j}: 1\le j \le k\}$.
The condition for the boundary map, $\partial_k\circ\partial_{k+1}=0$,
is then easily verified, and we have:
\begin{Theorem}
The map $\partial_k$ of (\ref{bmod}) defines a chain complex ${\cal M}_*$,
\begin{equation}
\label{chaincomplex}
0 \longrightarrow  {\cal M}_l~
{\mathop{\longrightarrow}^{\partial_l}} ~{\cal M}_{l-1}~
{\mathop{\longrightarrow}^{\partial_{l-1}}}~ \cdots~
{\mathop{\longrightarrow}^{\partial_2}}~ {\cal M}_1~
{\mathop{\longrightarrow}^{\partial_1}} ~{\cal M}_0~
{\longrightarrow}~0,
\end{equation}
where ${\cal M}_k$ is defined by (\ref{mk}).
\end{Theorem}
Since we have an explicit representation of the $k$-chains as (\ref{mk}) with
(\ref{wmod}), the integral homology $H_k(\hat H_{\Bbb R};{\Bbb Z})=
Ker \partial_k/Im \partial_{k+1}$ can be also computed.  However an explicit
formula
may be too complicated.

\section{Morse theory and Homology}
\renewcommand{\theequation}{4.\arabic{equation}}\setcounter{equation}{0}
\renewcommand{\thefigure}{4.\arabic{figure}}\setcounter{figure}{0}

\subsection{Morse theory}
The generalized Toda equation can be expressed as a gradient
flow on the adjoint orbit of $\frak g$ \cite{bloch:98}. Here we
 consider a Morse decomposition of the manifold $\hat H_{\Bbb R}$
based on the gradient structure of the Toda flow.
Each critical point of the Toda vector field can be parametrized by
a unique element of the Weyl group $W$.  Then we define the
{\it unstable} and {\it stable} Weyl subgroups for $a\in W$ as
\begin{eqnarray}
\label{unstble}
\left.
\begin{array}{ccc}
&  W^{u}(a)=W_{\Pi^u_a},\quad
& \Pi^u_a:=\{~{\alpha_i}\in \Pi~:~\ell(as_{\alpha_i})>\ell(a)~\}~ \\
&{}&{} \\
&  W^{s}(a)=W_{\Pi^s_a}, \quad
 &\Pi^s_a:= \{~{\alpha_i}\in \Pi~:~\ell(as_{\alpha_i})<\ell(a)~\}~
\end{array}
\right\}
\end{eqnarray}
where $\ell (a)$ denotes the length of $a$.
We will use the same notation for the unstable and stable manifolds
generated by the Toda vector field corresponding to the critical
point $a$. Thus, depending on the context,
$W^{u}(a)$, $W^{s}(a)$ denotes either a subgroup of $W$
or a submanifold of $\hat {H }_{\Bbb R }$.   We also introduce labels in
the Dynkin diagram
to characterize the critical point by assigning $``0"$ in the
$i$-th place in the diagram if $s_{\alpha_i}\in W^{s}(a)$,
and $``*"$ if $s_{\alpha_i}\in W^{u}(a)$.  For example,
in the case of ${\frak g} ={\frak{sl}}(6:{\Bbb R})$,
the element $[2143]\in W$ is labeled as $(0*0**)$, where
$[2143]:=s_{\alpha_2}s_{\alpha_1}s_{\alpha_4}s_{\alpha_3}$.
The $W^u([2143])$ is then the subgroup generated by $\{s_{\alpha_2},
s_{\alpha_4},
s_{\alpha_5}\}$ and is diffeomorphic to ${\Bbb R}^3$.
In terms of {\it handle body}, the critical point $[2143]$
is identified as the product $D^2\times D^3$ where $D^n$ is the
$n$-dimensional disc.

With this identification, we have the Morse decomposition
of the manifold ${\hat H}_{\Bbb R}$,
\begin{equation}
{\hat H}_{\Bbb R} ~= ~ \bigcup_{a\in W}~ W^{u}(a) ~.
\label{morse}
\end{equation}
The index of the critical point $a\in W$ is defined as
\begin{equation}
\label{index}
Ind(a) := dim W^{u}(a)=\Big|\Pi^u_a\Big|,
\end{equation}
which is also given by the number of $*$'s in the labeled
Dynkin diagram,  e.g. for ${\frak {sl}}(6:{\Bbb R})$,
$Ind([2143])=3$.

The Toda flow defines a (directed) graph which provides 1-dimensional
connections
among the critical points corresponding to the 1-dimensional flow
of a ${\frak {sl}}(2,{\Bbb R})$ subsystem.
We call the graph {\it Toda graph} and is defined by

\vskip 0.5cm
\noindent
{\bf Definition 4.1:} {\it Toda graph.}
A directed graph is called Toda graph if each vertex defined by
$\langle a\rangle=a^{-1}\langle e\rangle$ with $a\in W$ has the connections to
another vertex $\langle b_i\rangle$ by
\begin{equation*}
b_i=as_{\alpha_i}~, \quad \quad~~~~~~~~{\text{for}} \quad~~~~~~~~~~
i=1,\ldots,l.
\end{equation*}
The direction in the connection between two vertices $a$ and $b$ is defined by
\begin{equation*}
a \rightarrow b~, \quad \quad ~~~~~~~~~~ {\text{if}} \quad ~~~~~~~~~ \ell(a)
<\ell(b)~.
\end{equation*}

In order to construct a Morse complex, a vector field on the manifold
must satisfy the {\it Morse-Smale condition}, that is,  the  intersection
 between (the manifolds) $W^u(a)$ and $W^s(b)$ for the
critical points $a$ and $ b$ must be transversal. However the corresponding
intersections in the case
of the Toda lattice are, in general, NOT transversal. We have:

\vskip 0.5cm
\noindent
{\bf Definition 4.2:} {\it Transversal connection (algebraic version).}
A connection $a\to b$ is transversal if
\begin{enumerate}
\item  $\Big|\Pi^u_a\cap \Pi^s_b\Big|=Ind(a)-Ind(b)$,
\item  $\langle W^u(a), W^s(b)\rangle =W$,
\item  $\Big| aW^u(a)\cap bW^s(b)\Big|=\Big|W_{\Pi^u_a\cap\Pi^s_b}\Big|$ .
\end{enumerate}

\vskip 0.3cm

This definition is motivated by:

\begin{Theorem}
For the Toda lattice vector field, each CLOSURE, $\overline{W^{u}}(a)$, ($
\overline{W^{s}}(a)$) $a\in W$ of
the unstable (stable) manifold respectively, is smooth. Moreover each
smooth manifold  $\overline{ W^{u}}(a)$
($\overline{W^{s}}(a)$)
is orientable and produces a cycle if and only if the subgroup $W^{u}(a)$ (
$W^{s}(a)$) is abelian.
 A connection
$a\to b$ is transversal (Definition (4.2)) if and only if $\overline{
W^{u}}(a)$,
$\overline{ W^{s}}(b)$ intersect transversally. The intersection is
diffeomorphic to
a circle.

\end{Theorem}
The proof of the theorem can be obtained from the methods developed
in \cite{casian:99}, and the detail will be given elsewhere.

We call a  graph with vertices given by $W$ and oriented edges $a\to b$
 satisfying the (algebraic) transversality conditions above
 with $Ind(a)=Ind(b)+1$, a {\it Morse-Smale graph}, if
 in addition, i) there is a perturbation of the Toda lattice which is
Morse-Smale and has
  the same set of critical points ($W$),  ii)
 $a\to b$ only if the manifolds $W^{u}(a)$ and $W^{s}(b)$ for this new
vector field intersect
 transversally. A Morse-Smale vector field can be obtained
by a small smooth perturbation of the Toda lattice as in \cite{smale:61}.
 We have  confirmed that the conditions 1 through 3 in Definition
4.2 are  sufficient to determine uniquely a Morse-Smale graph in the
cases of ${\frak g}\cong A_l$ up to $l=3$. However this may not be true in
general.

We now define a boundary map on the chain ${\mathcal C}_*$ of the cells of
unstable Weyl groups $W^u(a)$, i.e.

\begin{equation}
\label{chain}
{\mathcal C}_*=\bigoplus_{k=0}^l {\mathcal C}_k,  ~~~~\quad \quad {\mathcal
C}_k=
\sum_{Ind(a)=k}{\Bbb Z}\langle a \rangle ~.
\end{equation}
where $\langle a\rangle$ is the cell corresponding to $W^u(a)$.
The chain ${\mathcal C}_k$ is the set of all cells $\langle a\rangle$
with the labeled Dynkin diagram having $k$ number of $*$'s.
The boundary map $\partial_k: {\mathcal C}_k \rightarrow {\mathcal C}_{k-1}$
is then defined by
\begin{equation}
{\partial_k}~ : ~ \langle a \rangle \longmapsto \partial_k \langle a \rangle
=\sum_{Ind(b)=k-1}[a;b]~ \langle b\rangle ~,
\label{dmap}
\end{equation}
where all the connections $a\to b$ are edges in the Morse-Smale graph,
and the {\it incidence} number $[a;b]$ is given by
\begin{equation}
\label{inc}
[a;b]=\left( 1+ (-1)^{\sigma [a;b]}\right)(-1)^{\ell(a^{-1}b)+i} ~
\end{equation}
with
\begin{eqnarray*}
 \sigma[a;b]=\Big|\{~j~:~\epsilon_j\to\epsilon'_j<0~,
~\alpha_j\in\Pi^u_a~\}\Big|,
\end{eqnarray*}
where the index $i$ is given by $\{s_{\alpha_i}\}= \Pi^u_a\cap\Pi^s_b$,
i.e. the Dynkin diagram corresponding to $\langle b\rangle$ has $0$
in the $i$-th place in addition to the $0$'s in $\langle a\rangle$.
The sign change $\epsilon_i\to \epsilon'_i$ under the connection
$a\to b$ is defined as follows:
\begin{equation}
a^{-1}b\cdot (\epsilon_1, \ldots,\epsilon_l)=(\epsilon'_1,\ldots,\epsilon'_l)
\label{signc}
\end{equation}
where the initial signs $\epsilon_i$'s are taken as
\begin{eqnarray*}
\epsilon_j=\left\{
\begin{array}{ccc}
+ ~~~ & \text{if} ~~ & {\alpha_j}\in \Pi^s_a\cup \Pi^u_b, \\
- ~~~ & \text{if} ~~ & {\alpha_j}\in \Pi^u_a\cap \Pi^s_b \\
\end{array}
\right.
\end{eqnarray*}
and $\epsilon'_j$ is defined as (\ref{ewj}).
With this definition, we can determine the change of orientations of the
hypersurfaces
for $b_{i}>0$ and $b_{i}<0$ parallel to the surface given by $b_{i}=0$
under the action of $x=a^{-1}b$.

\subsection{Example of $A_l={\frak{sl}}(l+1;{\Bbb R})$}

Let us first introduce the following elements of $W=S_{l+1}$, the symmetry
group
of order $l+1$;
\begin{equation}
\label{s}
s_{ij}:=s_{\alpha_i}\cdots s_{\alpha_j}=[i\cdots j]~,
\end{equation}
where the numbers $i\cdots j$ denotes the consective numbers between $i$
and $j$
for $1\le i,j\le l$.  Then for example some (but not all) of the connections
from the top
cell
$\langle e \rangle=(*\cdots *)$ to the cells labeled $(*\cdots * 0 *\cdots *)$
with $0$ in the $j$-th place are expressed by $s_{ij}:\langle e\rangle
\to \langle s_{ij}\rangle$ with
\begin{equation}
\label{sij}
\langle s_{ij}\rangle := s_{ij}^{-1}\cdot\langle e \rangle=s_{ji}\cdot
\langle e \rangle,  ~~~\quad~ \text{for} ~~~\quad 1\le i, j \le l ~.
\end{equation}
Note that the cell $\langle s_{ij}\rangle$ is isomorphic to $A_{j-1}\times
A_{l-j}$ as an unstable manifold $W^u(s_{ij})$ with index $l-1$ generated
by the Toda flows.  In particular, all the cells of $A_{l-1}$-type
are given by $\langle s_{i1} \rangle$ and $\langle s_{il}\rangle$
for $i=1,\ldots, l$.
We call the $k$-cells of $A_k$-type the {\it principal} part (of $k$-cells),
and those of the $A_{j_1}\times \cdots\times A_{j_n}$ with $j_1+\cdots +j_n=k,
n>1$ the {\it whisker} part.  The set of all the (principal) $k$-cells of
$A_k$-type is denoted as ${\cal A}_k$.
For example, the boundary of a principal $k$-cell
labeled by $(0\cdots 0*\cdots *0\cdots 0)$ with $l-k$ zeros are written in
the sum of the principal and the whisker parts of $(k-1)$-cells.
We then define a boundary map $\overset{\circ}{\partial}$ on the principal
$k$-cells
into the projection of the boundary map $\partial$ on
the principal parts of $(k-1)$-cells:
\begin{equation}
\label{bmap}
\overset{\circ}{\partial}_k~:~ {\cal A}_{k} \longrightarrow {\cal A}_{k-1}~.
\end{equation}
For the bounary of the top cell, we have
\begin{equation}
\label{topb}
\overset{\circ}{\partial}_l\langle e\rangle=\sum_{i=1}^l\Big(
[e;s_{i1}]\langle s_{i1}\rangle
 +[e;s_{il}]\langle s_{il}\rangle \Big) ~,
\end{equation}
where the incidence numbers are computed as
\begin{equation*}
\label{inctop}
[e;s_{i1}]=[e;s_{l-i+1,l}]=2(-1)^{i+1}(1-\delta_{il})~.
\end{equation*}
Thus the principal part of the boundary of the top cell
consists of $2(l-1)$ cells of $A_{l-1}$-type, and
the cells $\langle s_{1l}\rangle$ and $\langle s_{l1}\rangle$
are {\it not} in the part of the boundary.
We also note that  $\langle s_{1l}\rangle$ and
$\langle s_{l1}\rangle$ are only cells of $A_{l-1}$-type
separated from the others and invariant
under the subgroup generated by $W^u(s_{1l})$ for $\langle s_{1l}\rangle$
and $W^u(s_{l1})$ for $\langle s_{l1}\rangle$.
Then one can identify the cells which are not included in any parts
of the boundaries of $A_k$-type as in the following Proposition,
\begin{Proposition}
\label{prince}
All the cells which are free from the boundaries of cells with higher indices
are generated by the following commutative diagram starting
from $\langle a_{0,0} \rangle:=\langle e\rangle$,
\begin{equation}
\label{cda}
\begin{CD}
\langle a_{i,j}\rangle @>s_{l-j,i+1}>> \langle a_{i,j+1} \rangle\\
@Vs_{i+1,l-j}VV @VVs_{i+1,l-j-1}V \\
\langle a_{i+1,j}\rangle @>s_{l-j,i+2}>> \langle a_{i+1,j+1}\rangle
\end{CD}
\end{equation}
where $\langle a_{i,j}\rangle$ represents a unique cell labeled with
$(\overbrace{0\cdots 0}^i*\cdots *{\overbrace{0\cdots 0}^j}) \in
{\cal A}_{l-(i+j)}$.
\end{Proposition}
\begin{Proof}
The braid relation $[i\cdot i+1\cdot i]=[i+1\cdot i\cdot i+1]$ shows
the commutativity of the diagram. One can also show in a similar way as
in (\ref{topb}) that
the principal part of the boundary of the $\langle a_{i,j}\rangle$
consists of $2(l-(i+j)-1)$ cells and is given by
\begin{equation}
\overset{\circ}{\partial}_{l-(i+j)}\langle a_{i,j}\rangle=2
\sum_{k=1}^{l-(i+j)}(-1)^{k+1}\left(1-\delta_{k,l-(i+j)}\right)
 \Big(\langle b_{i+1,i+k}\rangle+\langle b_{l-j,l-j-k+1}\rangle\Big),
\label{kboundary}
\end{equation}
where $ \langle b_{i',j'}\rangle=s_{i',j'}\cdot\langle a_{i,j}\rangle$, and
the cells $\langle b_{i+1,l-j}\rangle=\langle a_{i+1,j}\rangle$
and $\langle b_{l-j,i+1}\rangle=\langle a_{i,j+1}\rangle$ do not appear.
\end{Proof}

The cells defined in Proposition \ref{prince} give the {\it seed} elements
of the ''principal graph" defined as the graph on the sets
${\cal A}_*:=\bigoplus_{k=1}^l{\cal A}_k$ where the connections
indicate the nonzero incidence numbers. Thus in the principal graph
there are $\displaystyle{{l(l+1) \over 2}}$ disconnected subgraphs,
each of which has a cell $\langle a_{i,j}\rangle$ as the highest dimensional
cell (the seed cell) in the subgraph with $dim\langle a_{i,j}\rangle
=l-(i+j)$.  Then one can show that the pair
$({\cal A}_*, \overset{\circ}{\partial}_*)$ forms a subchain complex, that is,
the boundary map satisfies $\overset{\circ}{\partial}_{k}\circ{\overset{\circ}
{\partial}_{k+1}}=0$.

\begin{figure}
\label{fig:3}
\epsfysize=2.8cm
\centerline{\epsffile{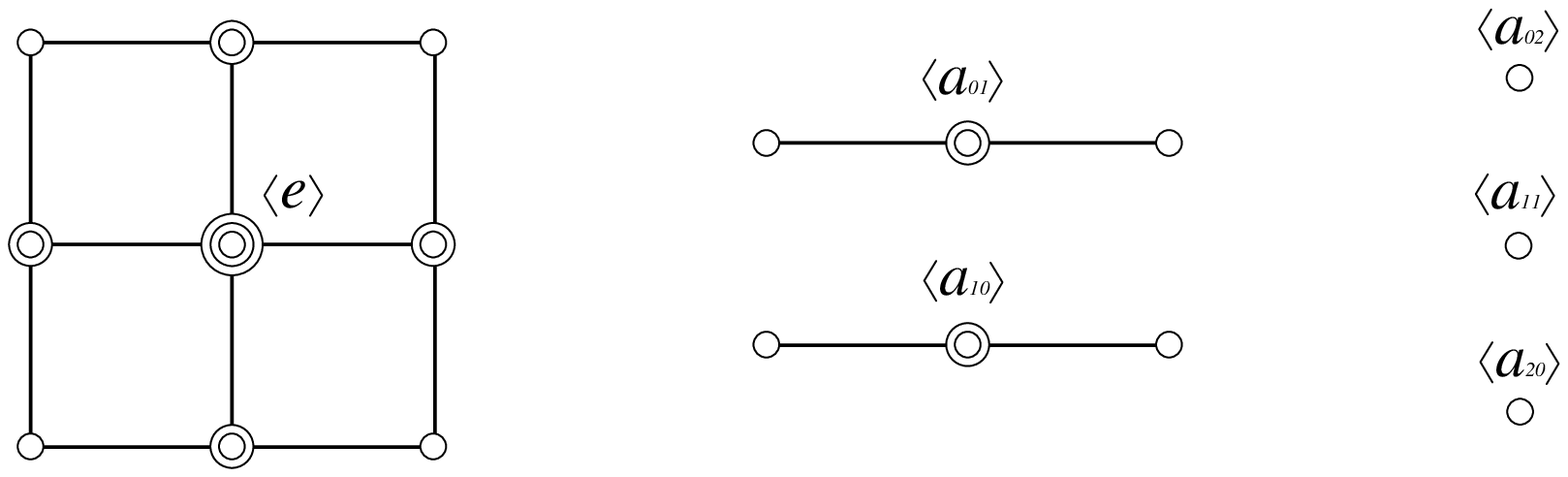}}
\caption{The principal graph for
$A_3={\frak{sl}}(4, {\Bbb R})$.}
\end{figure}

Figure 4.1 illustrates the example of $A_3={\frak {sl}}(4;{\Bbb R})$.
As in Fig.4.1, one can identify the cells in ${\cal A}_k$ as
$(k-1)$-dimensional cells in the graph consisting of $l(l+1)/2$ number
of disconnected {\it hypercubes}.  In each hypercube of dimension $k-1$, the
top cell
is represented by $\langle a_{i,j}\rangle$ with $k=l-(i+j)$, and
the vertices represent $A_1$-cycles.
For the case of $A_3$, we identify the seed cell $\langle e\rangle$ as the
face (square), the cells in ${\cal A}_2$ as the edges, and those in ${\cal
A}_1$
as the vertices. Then counting the numbers of those cells, we obtain
\begin{Theorem}
\label{h1complex}
The generating function (Poincar\'{e} polynomial) $P({\cal A}_*;q)$
 of the number of cells $|{\cal A}_k|$ is given by
\begin{equation}
\label{poincare}
P({\cal A}_*;q)=\sum_{k=1}^l|{\cal A}_k|q^{k-1}=\sum_{n=1}^l n(q+2)^{l-n}~.
\end{equation}
\end{Theorem}
\begin{Proof}
It is easy to see that the number of $k$-dimensional cells $a_k$ in the
$n$-dimensional hypercube is given by $\displaystyle{\binom{n}{k} 2^{n-k}}$
so that we have
\[
\sum_{k=0}^n |a_k|q^{k}=(q+2)^n~.
\]
From Proposition \ref{prince}, we have $n$ number of the $(l-n)$ dimensional
hypercubes in the principal graph.  This asserts the Theorem.
\end{Proof}

As a corollary of Theorem \ref{h1complex}, we obtain
\begin{Corollary}
\label{betti1}
The Betti number of $H_1({\hat H}_{\Bbb R},{\Bbb Z})$ is given by
\begin{equation}
\label{bett1}
b_1({\hat H}_{\Bbb R}):=rank(H_1)=P({\cal A}_*;-1)={l(l+1)\over 2}.
\end{equation}
\end{Corollary}
\begin{Proof}
Total number of $A_1$-cycles is given by the number of
vertices in the graph, i.e.
\[
|Z_1|=\sum_{n=0}^{l-1}2^n(l-n)=2^{l+1}-(l+2).
\]
From the graph, we can also find the number of boundaries,
that is, in each graph of $n$-dimensional hypercubes
there are $2^n-1$ boundaries, and we have $l-n$ disconnected
graph in this dimension.  Then we have
\[
|B_1|=\sum_{n=1}^{l-1}(l-n)\times (2^n-1),
\]
and obtain the Betti number $b_1=|Z_1|-|B_1|$ as stated.
\end{Proof}

Although we have a complete characterization of the cells in terms
of colored Dynkin diagrams (Section 3), we have not obtained
explicitly higher homology.  It is however
natural to consider the following conjecture on the Betti
numbers $b_k$ as the alternative sums of the numbers of
{\it whiskers}:
\begin{eqnarray}
\label{bettik}
b_k=
\left\{
\begin{array}{ll}
\displaystyle{\sum_{n=k}^{l-k+1} \Big|{\cal A}_n^{(k)}\Big|(-1)^{n-k}},
\quad\quad &{\mbox{for}}\quad 1\le k \le {l+1\over 2}, \\
\quad\quad \quad 0 ~, \quad\quad & {\mbox{for}} \quad  k>{l+1\over 2},
\end{array}
\right.
\end{eqnarray}
where $|{\cal A}_n^{(k)}|$ is the number of whiskers defined by
\[
\Big|{\cal A}_n^{(k)}\Big|:=
\sum_{n_1+\cdots+n_k=n \atop 1\le n_1\le \cdots\le n_k}
\Big|{\cal A}_{n_1}\times \cdots \times {\cal A}_{n_k}\Big|.
\]
Note here that all the $k$-cycles are given by the products of $A_1$-cycles,
i.e. $|Z_k|=|{\cal A}_k^{(k)}|$.  The conjecture is confirmed
for the cases of ${\frak g}\cong A_l$ up to $l=3$.

\section{Final remark}

In this paper, we have studied the topology of the isospectral manifolds
 associated with the compactified level variety of the generalized Toda
(Kostant-Toda) lattices on real split semisimple Lie algebras.  The details
of the decomposition based on the
colored Dynkin diagrams can be found in our recent paper \cite{casian:99},
and the proofs of the results stated in the sections 2 and 3 can be also
found in this paper.

As a final remark, we would like to mention a possible extension
of the present study for the
full Kostant-Toda lattices which are recently shown to be integrable
in \cite{gekhtman:99}. Our methods may then shed some light on the structure
of the real full flag manifold.

\bibliographystyle{amsplain}

\end{document}